\newcount\secno
\newcount\prmno
\newif\ifnotfound
\newif\iffound

\def\namedef#1{\expandafter\def\csname #1\endcsname}
\def\nameuse#1{\csname #1\endcsname}

\long\def\ifundefined#1#2#3{\expandafter\ifx\csname
  #1\endcsname\relax#2\else#3\fi}
\def\hwrite#1#2{{\let\the=0\edef\next{\write#1{#2}}\next}}

\toksdef\ta=0 \toksdef\tb=2
\long\def\leftappenditem#1\to#2{\ta={\\{#1}}\tb=\expandafter{#2}%
                                \edef#2{\the\ta\the\tb}}
\long\def\rightappenditem#1\to#2{\ta={\\{#1}}\tb=\expandafter{#2}%
                                \edef#2{\the\tb\the\ta}}

\def\lop#1\to#2{\expandafter\lopoff#1\lopoff#1#2}
\long\def\lopoff\\#1#2\lopoff#3#4{\def#4{#1}\def#3{#2}}

\def\ismember#1\of#2{\foundfalse{\let\given=#1%
    \def\\##1{\def\next{##1}%
    \ifx\next\given{\global\foundtrue}\fi}#2}}

\def\section#1{\vskip1cm
               \global\def\currenvir{section}
               \global\advance\secno by1\global\prmno=0
               {\bf \number\secno. {#1}}
               \smallskip}

\def\subsection{\global\def\currenvir{subsection}
                \global\advance\prmno by1
                \ind{ (\number\secno.\number\prmno) }}
\def\subsec{\global\def\currenvir{subsection}
                \global\advance\prmno by1
                { (\number\secno.\number\prmno)\ }}

\def\proclaim#1{\global\advance\prmno by 1
                {\bf #1 \the\secno.\the\prmno$.-$ }}

\long\def\th#1 \enonce#2\endth{%
   \medbreak\proclaim{#1}{\it #2}\global\def\currenvir{th}\smallskip}

\def\rem#1{\global\advance\prmno by 1
{\it #1}\kern5pt\the\secno.\the\prmno$.-$ }


\def\isinlabellist#1\of#2{\notfoundtrue%
   {\def\given{#1}%
    \def\\##1{\def\next{##1}%
    \lop\next\to\za\lop\next\to\zb%
    \ifx\za\given{\zb\global\notfoundfalse}\fi}#2}%
    \ifnotfound{\immediate\write16%
                 {Warning - [Page \the\pageno] {#1} No reference found}}%
                \fi}%
\def\ref#1{\ifx\labellist\empty{\immediate\write16
                 {Warning - No references found at all.}}
               \else{\isinlabellist{#1}\of\labellist}\fi}

\def\newlabel#1#2{\rightappenditem{\\{#1}\\{#2}}\to\labellist}
\def\labellist{}

\def\label#1{\relax}
 \def\openall{\openout\lbl=\jobname.lbl}

\newread\testfile
\def\lookatfile#1{\openin\testfile=\jobname.#1
    \ifeof\testfile{\immediate\openout\nameuse{#1}\jobname.#1
                    \write\nameuse{#1}{}
                    \immediate\closeout\nameuse{#1}}\fi%
    \immediate\closein\testfile}%

\def\begin{  \newlabel{e}{1.1}
\newlabel{free}{2.1}
\newlabel{g}{2.2}
\newlabel{conc}{2.3}
\newlabel{embed}{3.1}
\newlabel{elag}{3.3}
\newlabel{An}{3.5}
\newlabel{homeo}{3.7}
\newlabel{prod}{3.8}
\newlabel{C*}{4.1}
\newlabel{\Gamma }{4.2}
\newlabel{pq}{4.3}
\newlabel{rat}{4.5}}


\magnification 1250
\pretolerance=500 \tolerance=1000  \brokenpenalty=5000
\mathcode`A="7041 \mathcode`B="7042 \mathcode`C="7043
\mathcode`D="7044 \mathcode`E="7045 \mathcode`F="7046
\mathcode`G="7047 \mathcode`H="7048 \mathcode`I="7049
\mathcode`J="704A \mathcode`K="704B \mathcode`L="704C
\mathcode`M="704D \mathcode`N="704E \mathcode`O="704F
\mathcode`P="7050 \mathcode`Q="7051 \mathcode`R="7052
\mathcode`S="7053 \mathcode`T="7054 \mathcode`U="7055
\mathcode`V="7056 \mathcode`W="7057 \mathcode`X="7058
\mathcode`Y="7059 \mathcode`Z="705A
\def\spacedmath#1{\def\packedmath##1${\bgroup\mathsurround =0pt##1\egroup$}
\mathsurround#1
\everymath={\packedmath}\everydisplay={\mathsurround=0pt}}
\def\nospacedmath{\mathsurround=0pt
\everymath={}\everydisplay={} } \spacedmath{2pt}

\def\phfl#1#2{\normalbaselines{\baselineskip=0pt
\lineskip=10truept\lineskiplimit=1truept}\nospacedmath\smash {\mathop{\hbox to
8truemm{\rightarrowfill}}
\limits^{\scriptstyle#1}_{\scriptstyle#2}}}
\def\hfl#1#2{\normalbaselines{\baselineskip=0truept
\lineskip=10truept\lineskiplimit=1truept}\nospacedmath\smash{\mathop{\hbox to
12truemm{\rightarrowfill}}\limits^{\scriptstyle#1}_{\scriptstyle#2}}}
\def\diagram#1{\def\normalbaselines{\baselineskip=0truept
\lineskip=10truept\lineskiplimit=1truept}   \matrix{#1}}
\def\vfl#1#2{\llap{$\scriptstyle#1$}\left\downarrow\vbox to
6truemm{}\right.\rlap{$\scriptstyle#2$}}
\font\eightrm=cmr8\font\sixrm=cmr6
\def\note#1#2{\footnote{\parindent
0.4cm$^#1$}{\vtop{\eightrm\baselineskip12pt\hsize15.5truecm\noindent #2}}
\parindent 0cm}

\def\iso{\mathrel{\mathop{\kern 0pt\longrightarrow }\limits^{\sim}}}
\def\union_#1^#2{\mathrel{\mathop{\kern0pt\scriptscriptstyle\bigcup}
\limits_{#1}^{#2}}}
\def\sdir_#1^#2{\mathrel{\mathop{\kern0pt\oplus}\limits_{#1}^{#2}}}
\def\pprod_#1^#2{\raise 2pt
\hbox{$\mathrel{\scriptstyle\mathop{\kern0pt\prod}\limits_{#1}^{#2}}$}}
\def\cprod_#1^#2{\raise 2pt
\hbox{$\mathrel{\scriptscriptstyle\mathop{\kern0pt\coprod}\limits_{#1}^{#2}}$}}

\def\pc#1{\tenrm#1\sevenrm}
\def\up#1{\raise 1ex\hbox{\smallf@nt#1}}
\def\tx{\kern-1.5pt -}
\def\cqfd{\kern 1truemm\unskip\penalty 500\vrule height 4pt depth 0pt width
4pt\medbreak} 
\def\virg{\raise
.4ex\hbox{,}}
\def\decale#1{\smallbreak\hskip 28pt\llap{#1}\kern 5pt}
\def\no{n\up{o}\kern 2pt}
\def\ind{\par\hskip 1truecm\relax}
\def\indp{\par\hskip 0.5cm\relax}
\def\moins{\mathrel{\hbox{\vrule height 3pt depth -2pt width 6pt}}}
\def\rond{\kern 1pt{\scriptstyle\circ}\kern 1pt}
\def\iso{\mathrel{\mathop{\kern 0pt\longrightarrow }\limits^{\sim}}}

\def\Hom{\mathop{\rm Hom}\nolimits}
\def\Aut{\mathop{\rm Aut}\nolimits}

\def\Pic{\mathop{\rm Pic}\nolimits}

\def\dim{\mathop{\rm dim}\nolimits}
\def\Card{\mathop{\rm Card}\nolimits}

\def\jd{\bar{\cal J}^d{\cal C}}
\def\jg{\bar{\cal J}^g{\cal C}}
\font\teneufm=eufm10
\newfam\gothfam
\textfont\gothfam=\teneufm
\def\goth{\fam\gothfam}
\vsize = 25.3truecm
\hsize = 16truecm
\voffset = -.5truecm
\parindent=0cm
\baselineskip15pt
\begin
\centerline{\bf Counting rational curves on K3 surfaces}
\smallskip
\smallskip \centerline{Arnaud {\pc BEAUVILLE\note{1}{Partially supported by the
European HCM project ``Algebraic Geometry in Europe" (AGE).}}}
\smallskip
\centerline{Version 0}
\vskip1.2cm

{\bf Introduction}
\smallskip
\ind The aim of these notes is to explain the remarkable formula found by Yau
and Zaslow [Y-Z] to express the number of rational curves on a K3 surface.
Projective K3 surfaces fall into countably many
families $({\cal F}_g)_{g\ge1}$; a surface in ${\cal F}_g$
admits a
$g$\tx dimensional linear system of curves of genus $g$. A na\"{\i}ve count of
constants suggests that such a system will contain a positive number, say
$n(g)$, of rational (highly singular) curves. The formula is
$$\sum_{g\ge0}n(g)q^g= {q\over \Delta (q)}\ ,$$
where $\Delta (q)=q\prod_{n\ge1}(1-q^n)^{24}$ is the well-known modular form
of weight $12$, and we put by convention $n(0)=1$.
\ind To explain the idea in a nutshell, take the case
$g=1$. We are thus looking at  K3 surfaces with an elliptic
fibration $f:S\rightarrow {\bf P}^1$, and we are asking for the number of
singular fibres. The (topological) Euler-Poincar\'e characteristic of a fibre
$C_t$ is $0$ if $C_t$ is smooth, $1$ if it is a  rational curve with one
node, $2$ if it has a cusp, etc. From the standard properties of the
Euler-Poincar\'e  characteristic, we get
$\displaystyle e(S)=\sum_t e(C_t)$; hence $n(1)=e(S)=24$, and this number
counts
nodal rational curves with multiplicity $1$, cuspidal rational curves with
multiplicity $2$, etc. \ind The idea of Yau and Zaslow is to generalize this
approach to any genus. Let $S$ be a K3 surface with a $g$\tx dimensional linear
system $\Pi$ of curves of genus $g$. The role of $f$ will be played by the
morphism $\bar{\cal J}{\cal C} \rightarrow \Pi$ whose fibre over a point
$t\in\Pi$
is the compactified Jacobian $\bar JC_t$. To apply the same method, we would
like to prove the following facts:

\indp 1) The Euler-Poincar\'e characteristic $e(\bar{\cal J}{\cal C})$ is the
coefficient of $q^g$ in the Taylor expansion of $q/\Delta(q)$.
\indp 2) $e(\bar JC_t)=0$  if $C_t$ is not rational.

\indp 3) $e(\bar JC_t)=1$ if $C_t$ is a rational curve with nodes as only
singularities. Moreover $e(\bar JC_t)$ is positive when $C_t$ is rational, and
can
be computed in terms of the singularities of $C_t$.

\indp 4) For a generic K3 surface $S$ in ${\cal F}_g$, all rational curves in
$\Pi $
are nodal.

 \ind The first statement is proved in \S 1, by comparing with
the Euler-Poincar\'e characteristic of the Hilbert scheme $S^{[g]}$ which has
been
computed by G\"ottsche. The assertion 2) is proved in \S 2. The situation for
3)
is less satisfactory: though I can express $e(\bar JC)$, for a rational curve
$C$,
in terms of a local invariant of the singularities of $C$, and compute this
local
invariant in a number of cases, at this moment I am
 not able to prove that it is always positive.
 Finally 4) seems to be still open, despite recent progress by Xi Chen.
\ind The outcome (see Cor.\ \ref{conc}) is that the coefficient of $q^g$ in
$q/\Delta (q)$ counts the rational curves in $\Pi$ with a certain multiplicity,
which is $1$ for a nodal curve and can be computed explicitely in many cases;
the
two missing points are the positivity of this multiplicity, and the fact that
only
nodal curves occur on a generic surface in ${\cal F}_g$.

 \section{The compactified relative Jacobian }
\subsection\label{e}  Let $X$ be a complex variety; we denote by $e(X)$ its
Euler-Poincar\'e characteristic, defined by $e(X)=\sum_p(-1)^p\dim_{\bf
Q}H_c^p(X,{\bf Q})$. Recall that this invariant is additive, i.e.\ satisfies
 $e(X)=e(U)+e(X\moins
U)$ whenever  $U$ is an open subset of $X$.
\subsection We consider a projective K3 surface $S$
with a complete linear system $(C_t)_{t\in \Pi}$ of curves of genus $g\ge1$ (so
$\Pi $ is a projective space of dimension $g$). We will  assume that {\it all
the
curves $C_t$ are integral} (i.e.\ irreducible and reduced). This is a
simplifying
assumption, which can probably be removed at the cost of various technical
complications. It is of course satisfied if the class of $C_t$ generates
$\Pic(S)$.  \ind  Let ${\cal C}\rightarrow \Pi $ be the morphism with fibre
$C_t$
over $t\in\Pi$. For each integer $d\in {\bf Z}$, we denote by $\bar{\cal
J}{\cal C}=\cprod_{d\in{\bf Z}}^{}\bar{\cal J}^d{\cal C}$ the compactified
Picard scheme of this family. $\jd$ is a projective
variety of dimension $2g$, which parameterizes pairs $(C_t,{\cal L})$ where
$t\in\Pi $ and
 ${\cal L}$ is a torsion
free, rank $1$ coherent sheaf on $C_t$ of degree $d$ (i.e.\ with $\chi ({\cal
L})=d+1-g$). According to Mukai ([M], example 0.5), $\jd$ can be viewed as a
connected component of the moduli space of simple sheaves on $S$, and
therefore is smooth, and admits a (holomorphic) symplectic
structure.

\ind The simplest symplectic varieties associated to the K3 surface $S$ are
the Hilbert schemes $S^{[d]}$, which parameterize finite subschemes of
length $d$ of $S$.  The birational comparison of the symplectic varieties
$\jd$,
for various values of
$d$, with $S^{[g]}$ is an interesting problem, about which not much seems to be
known. There is one easy case:
\th Proposition
\enonce The compactified Jacobian $\jg$ is birationally isomorphic
to~$S^{[g]}$.
\endth
{\it Proof}: Let $U$ be the open subset of $\jg$ consisting of pairs
$(C_t,L)$ where $L$ is invertible and $\dim H^0(C_t,L)=1$. To such a pair
corresponds a unique effective Cartier divisor $D$ on $C_t$ of degree $g$,
which
can be viewed as a length $g$  subscheme  of $S$; since
 $\dim H^0(C_t,{\cal O}_{C_t}(D))=1$ it is contained in a unique curve
of $\Pi$, namely $C_t$. This provides an isomorphism betwen $U$ and the open
subset $V$ of $S^{[g]}$ parameterizing  finite subschemes of $S$ contained in a
unique curve of $\Pi$ and defining a Cartier divisor in this curve. \cqfd
\smallskip
\th Corollary
\enonce   Write $\displaystyle {q\over \Delta(q)}=\sum_{g\ge0}e(g)\,q^g$. Then
$e(\jg)=e(g)$.
\endth
{\it Proof}: We can either use a recent result of Batyrev and Kontsevich [?]
saying that two birationnally equivalent projective Calabi-Yau manifolds have
the
same Betti numbers, or a more precise result of Huybrechts [H]:
 two birationally equivalent
projective  symplectic manifolds which are isomorphic in codimension $2$ are
diffeomorphic (note that the open subsets $U$ and $V$ appearing in the above
proof
have complements of codimension $\ge2$). It remains to apply G\"ottsche's
formula
$e(S^{[g]})=e(g)$ [G]. \cqfd

\section{The compactified Jacobian of a non-rational curve}\smallskip
\ind Let $C$ be an integral curve.
By a {\it rank $1$ sheaf} on $C$ I will mean
a torsion free, rank $1$ coherent sheaf.
The  rank $1$  sheaves ${\cal L}$ on
$C$ of degree $d$ are parameterized by the
compactified Jacobian
$\bar J^dC$. If $L$ is an invertible sheaf of degree $d$ on $C$, the map
${\cal L}\mapsto {\cal L}\otimes L$ is an isomorphism of $\bar J C$ onto $\bar
J^dC$, so we can restrict our study to degree $0$ sheaves.
\ind Let ${\cal L}\in \bar JC$; the endomorphism ring of ${\cal L}$ is an
${\cal O}_C$\tx subalgebra of the sheaf of rational functions on $C$. It is
finitely generated as a ${\cal O}_C$\tx module, hence contained in ${\cal
O}_{\widetilde{C}}$. It is thus of the form ${\cal O}_{C'}$, where
$f:C'\rightarrow C$ is some partial normalization of $C$. The sheaf ${\cal L}$
is
a ${\cal O}_{C'}$\tx module, which amounts to say that it is the direct image
of a
rank $1$ sheaf ${\cal L}'$ on $C'$.

\th Lemma
\enonce Let $L\in JC$. Then ${\cal L}\otimes L$ is isomorphic to ${\cal L}$ if
and
only if $f^*L$ is trivial.
\endth\label{free}
{\it Proof}: The sheaf ${\cal L}\otimes L$ is isomorphic to $f_*({\cal
L}'\otimes f^*L)$, hence to ${\cal L}$ if $f^*L$ is trivial.
On the other hand we have
$${\cal H}om_{{\cal O}_C}({\cal L},{\cal L}\otimes L)\cong {\cal E}nd_{{\cal
O}_C}({\cal L})\otimes_{{\cal O}_C}L\cong f_*{\cal O}_{C'}\otimes L\cong
f_*f^*L\ ,$$
so if $f^*L$ is non-trivial, the space $\Hom({\cal L},{\cal L}\otimes L)$ is
zero, and ${\cal L}\otimes L$ cannot be isomorphic to ${\cal L}$. \cqfd
\smallskip

\th Proposition
\enonce Let $C$ be an integral curve whose normalization $\widetilde{C}$ has
genus $\ge1$. Then $e(\bar J^d C)=0$.
\endth\label{g}
{\it Proof}: We have an exact sequence
$$0\rightarrow G\longrightarrow JC\longrightarrow J\widetilde{C}\rightarrow
0\ ,$$
where $G$ is a product of additive and multiplicative groups. In particular,
$G$ is a divisible group, hence this exact sequence splits as a sequence of
abelian groups. For each integer $n$, we can therefore find a subgroup of order
$n$
 in
$JC$ which maps injectively into  $J\widetilde{C}$.
By Lemma
\ref{free}, this group acts freely on $\bar JC$, which implies that $n$ divides
$e(\bar JC)$; since this holds for any $n$ the Proposition follows. \cqfd
\smallskip
\th Corollary
\enonce Write $\displaystyle {q\over \Delta(q)}=\sum_{g\ge0}e(g)\,q^g$; let
$\Pi_{\rm rat}\i\Pi$ be the (finite) subset of  rational curves. Then
$\displaystyle e(g)=\sum_{t\in\Pi_{\rm rat}}e(\bar JC_t)$.
\endth\label{conc}
{\it Proof}: We first make a general observation: let $f:X\rightarrow Y$ be a
surjective morphism of complex algebraic varieties whose fibres have Euler
characteristic $0$; then $e(X)=0$. This is well known (and easy) if $f$ is a
locally trivial fibration; the general case follows using (\ref{e}), because
there
exists a stratification of $Y$ such that $f$ is locally trivial above each
stratum
[V]. \ind The set
$\Pi_{\rm rat}$ is finite because otherwise it would contain a curve, so $S$
would be
ruled. Consider the morphism $p:\jg\rightarrow \Pi$ above $\Pi\moins\Pi_{\rm
rat}$;
by the above remark, we have $e(p^{-1} (\Pi\moins\Pi_{\rm rat}))=0$, hence the
result using again (\ref{e}). \cqfd \medskip
\ind In other words, $e(g)$ counts the
number of rational curves with multiplicity, the multiplicity of a curve $C$
being
$e(\bar JC)$. In the next two sections we will try to show that this is indeed
a
reasonable notion of multiplicity (with only partial success, as explained in
the
introduction).

\section{The compactified Jacobian of a rational curve}
\th Lemma
\enonce Let $f:C'\rightarrow C$ be a partial normalization of $C$. The morphism
$f_*:\bar JC'\rightarrow \bar JC$ is a closed embedding.
\endth\label{embed}
{\it Proof}: Let ${\cal L},{\cal M}$ be two rank $1$ sheaves on $C'$.
We claim that any ${\cal O}_C$\tx homo\-morphism $u:f_*{\cal L} \rightarrow
f_*{\cal M}$ is actually $f_*{\cal O}_{C'}$\tx linear. Let $U$ be a Zariski
open
subset of $C$, $\varphi\in \Gamma (U,f_*{\cal O}_{C'}) $, $s\in \Gamma
(U,f_*{\cal
L})$; the rational function $\varphi $ can be written as $a/b$, with $a,b\in
\Gamma (U,{\cal O}_{C})$ and $b\not=0$. Then the element $u(\varphi s)-\varphi
u(s)$ of $\Gamma (U,f_*{\cal M})$ is killed by $b$, hence is zero since
$f_*{\cal
M}$ is torsion-free.
\ind Therefore if $f_*{\cal L}$ and $f_*{\cal M}$ are isomorphic as ${\cal
O}_C$\tx modules, they are also isomorphic as $f_*{\cal O}_{C'}$\tx modules,
which
means that ${\cal L}$ and ${\cal M}$ are isomorphic: this proves the
injectivity
of $f_*$ (which would be enough for our purpose). Now if $S$ is any base
scheme,
the same  argument applies to sheaves ${\cal L}$, ${\cal M}$  on $C\times S$,
flat
over $S$, whose restrictions to each fibre $C\times\{s\}$ are torsion free rank
$1$ (observe that a local section $b$ of ${\cal O}_C$ is ${\cal M}$\tx
regular because it is on each fibre, and ${\cal M}$ is flat over $S$). This
proves that $f_*$ is a monomorphism; since it is proper, it is a closed
embedding.
\cqfd

\medskip

\subsection Recall that the curve $C$ is said to
be {\it unibranch} if its normalization $\widetilde{C}\rightarrow C$ is a
homeomorphism. Any curve $C$ admits a unibranch partial normalization $\check
\pi:
\check C\rightarrow C$ which is minimal, in the sense   that any unibranch
partial normalization $C'\rightarrow C$ factors through $\check \pi$. To see
this, let
${\cal C}$ be the conductor of $C$, and let $\widetilde{\Sigma }$ be the
inverse
image in $\widetilde{C}$ of the singular locus $\Sigma \in C$. The
finite-dimensional $k$\tx algebra $A:={\cal O}_{\widetilde{C}}/{\cal C}$ is a
product of local  rings $(A_{x})_{x\in\widetilde{\Sigma }}$;  let
$(e_x)_{x\in\widetilde{\Sigma }}$ be the corresponding idempotent elements of
$A$.
A sheaf of algebras ${\cal O}_{C'}$ with ${\cal O}_{C}\i{\cal O}_{C'}\i{\cal
O}_{\widetilde{C}}$ is unibranch if and only if ${\cal O}_{C'}/{\cal C}$
contains
each $e_x$, or equivalently  ${\cal O}_{C'}$ contains the  classes $e_x+{\cal
C}$
for each $x\in\widetilde{\Sigma }$; clearly there is a smallest such algebra,
namely the algebra  ${\cal O}_{\check C}$  generated by ${\cal O}_C$ and the
classes $e_x+{\cal C}$. The completion of the local ring of $\check C$ at a
point
$y$ is the image of $\widehat{\cal O}_{C,\check \pi(y)}$ in $\widehat{\cal
O}_{\widetilde{C},y}$.

\smallskip

\th Proposition
\enonce With the above notation, $e(\bar JC)=e(\bar J\check C)$.
\endth \label{elag}
{\it Proof}: In view of Prop.\ \ref{g}, we may suppose that $\widetilde{C}$ is
rational.  As
before we denote by
$\Sigma $ the singular locus of $C$, and by $\widetilde{\Sigma }$ its inverse
image
in $\check C$.
The cohomology exact sequence associated to the short exact sequence
$$1\rightarrow {\cal O}_C^*\longrightarrow{\cal O}_{\widetilde{C}}^*
\longrightarrow {\cal O}_{\widetilde{C}}^*/{\cal O}_C^*\rightarrow 1$$
provides a bijective homomorphism (actually an isomorphism of algebraic groups)
${\cal O}_{\widetilde{C}}^*/{\cal O}_C^*\iso JC$.

 \ind  The evaluation maps
${\cal O}_{\widetilde{C}}^*\rightarrow ({\bf C}^*)^{\widetilde{\Sigma}} $ and
${\cal O}_C^*\rightarrow ({\bf C}^*)^{\Sigma} $ give rise to a surjective
homomorphism ${\cal O}_{\widetilde{C}}^*/{\cal O}_C^*\rightarrow ({\bf
C}^*)^{\widetilde{\Sigma}}/({\bf C}^*)^{\Sigma}$; its kernel is unipotent,
i.e.\
isomorphic to a vector space. If $n$ is any integer $\ge\Card(\widetilde{\Sigma
})$, it follows that we can find a section $\varphi$ of ${\cal
O}_{\widetilde{C}}^*$ in a neighborhood of $\widetilde{\Sigma} $ such that the
numbers $\varphi(\tilde x)$ for $\tilde x\in\widetilde{\Sigma}$ are all
distinct,
but $\varphi ^n$ belongs to ${\cal O}_C$. Let $L$ be the line bundle on $JC$
associated to the class of $\varphi $ in ${\cal O}_{\widetilde{C}}^*/{\cal
O}_C^*$.

\ind Let $U$ be the complement of $\check \pi_*(\bar J\check C)$ in $\bar JC$;
according to  \ref{e} and Lemma \ref{embed}, our assertion is equivalent to
$e(U)=0$. We claim that the line bundle $L$ acts freely on $U$; since the order
of
$L$ in $JC$ is finite and arbitrary large, this will finish the proof. Let
${\cal
L}\in U$, and let $C'$ be the partial normalization of $C$ such that ${\cal
E}nd({\cal L})={\cal O}_{C'}$; by definition of $U$, $C'$ is not unibranch,
hence
there are two points of $\widetilde{\Sigma}$ mapping to the same point of $C'$;
this implies that the function $\varphi$ does not belong to ${\cal O}_{C'}^*$.
{}From the commutative diagram
$$\diagram{\widetilde{\cal O}_{\widetilde{C}}^*/{\cal O}_C^* &\hfl{\raise
-3mm\hbox{$\sim$}}{} &JC\cr \vfl{}{}& &\vfl{}{}\cr {\cal
O}_{\widetilde{C}}^*/{\cal O}_{C'}^* &\hfl{\sim}{} & JC' }$$
we conclude that the pull back  of $L$ to $JC'$ is non-trivial; by Lemma
\ref{free} this implies that ${\cal L}\otimes L$ is not isomorphic to ${\cal
L}$.
\cqfd

\th Corollary
\enonce For a rational nodal curve $C$, we have $e(\bar JC)=1$. \cqfd
\endth
\medskip
\rem{Remark}\label{An} Consider a rational curve $C$ whose  singularities are
all
of type $A_{2l-1}$, i.e.\ locally defined by an equation $u^2-v^{2l}=0$.
Locally around such a singularity, the curve $C$ is the union of two smooth
branches with a high order contact, so by \ref{elag} $e(\bar JC)$ is equal to
$1$.  The fact that some highly singular curves count with multiplicity one
looks
rather surprising. The case $g=2$ provides a (modest) confirmation: the surface
$S$ is a double covering of ${\bf P}^2$ branched along a sextic curve $B $; the
curves $C_t$  are the  inverse images of the lines in ${\bf P}^2$, and they
become
rational when the line is bitangent to $B$. We get an $A_3$\tx singularity when
the line has a contact of order $4$; thus our assertion in this case follows
from
the (certainly classical) fact that a line with a fourth order contact counts
as a
simple bitangent.

\medskip

\subsection Prop.\ \ref{elag} reduces the computation of the invariant $e(\bar
J
C)$ to the case of a unibranch (rational) curve. To understand this invariant
we
will use a  construction of Rego ([R], see also [G-P]).  For each $x\in C$, we
put
$\delta_x=\dim {\cal O}_{\widetilde{C},x}/{\cal O}_{C,x}$ and we denote by
${\cal
C}$ the ideal ${\cal O}_{\widetilde{C}}(-\sum_x (2\delta_x)\, x)$; it is
contained
in the conductor of $C$ (but the inclusion is strict unless $C$ is Gorenstein).
\ind For $x\in C$, we denote by $A_x$
and
$\widetilde{A}_x$ the finite dimensional algebras
${\cal O}_{C,x}/{\cal C}_x$ and ${\cal O}_{\widetilde{C},x}/{\cal C}_x$. Let
 ${\bf G}(\delta_x,\widetilde{A}_x)$  be the Grassmannian of codimension
$\delta_x$ subspaces of
$\widetilde{A}_x$, and  ${\bf G}_x$  the closed
subvariety of
 ${\bf G}(\delta_x,\widetilde{A}_x)$ consisting of sub-$\!A_x$\tx modules.
 We can also view  ${\bf
G}_x$ as parameterizing the  sub-$\!{\cal O}_{C,x}$\tx modules
${\cal L}_x$ of codimension $\delta_x$ in
${\cal O}_{\widetilde{C},x}$, because any such sub-module contains ${\cal C}_x$
([G-P], lemma 1.1 (iv)).
Since
${\cal O}_{\widetilde{C}}/{\cal C}$ is a skyscraper sheaf with fibre
$\widetilde{A}_x$ at
$x$,  the product $\pprod_{x\in\Sigma }^{}{\bf
G}_x$ parameterizes sub-$\!{\cal O}_C$\tx modules ${\cal L}\i
{\cal O}_{\widetilde{C}}$ such that $\dim {\cal O}_{\widetilde{C},x}/{\cal
L}_x=\delta_x$ for all $x$.
This implies $\chi ({\cal O}_{\widetilde{C}}/{\cal L})=\sum_x \delta_x
=\chi ({\cal O}_{\widetilde{C}}/{\cal O}_C)$, hence ${\cal L}\in \bar JC$.
We have thus defined a morphism $e:\pprod_{x\in\Sigma }^{}{\bf
G}_x\rightarrow \bar JC$.
\th Proposition
\enonce The map $e$ is a homeomorphism.
\endth\label{homeo}
\ind Note that $e$ is not an isomorphism, already when $C$ is a rational curve
with one ordinary cusp $s$:  the Grassmannian ${\bf
G}_s$ is isomorphic to ${\bf P}^1$, while $\bar JC$ is
isomorphic to $C$.
\ind  Since we are dealing with compact varieties, it suffices to prove
that
$e$ is bijective.\smallskip
 {\it Injectivity}: Let ${\cal L}$, ${\cal M}$ be two sub-$\!{\cal O}_C$\tx
modules
of ${\cal O}_{\widetilde{C}}$ containing ${\cal C}$. If ${\cal L}$ and ${\cal
M}$
give the same element in $\bar JC$, there exists a rational function $\varphi$
on
$\widetilde{C}$ such that ${\cal M}=\varphi{\cal L}$. But the equalities $\dim
{\cal O}_{\widetilde{C},x}/{\cal M}_x = \dim {\cal O}_{\widetilde{C},x}/{\cal
L}_x=\dim \varphi_x{\cal O}_{\widetilde{C},x}/{\cal M}_x$ imply $\varphi_x
{\cal
O}_{\widetilde{C},x}={\cal O}_{\widetilde{C},x}$ for all $x$, which means that
$\varphi$ is constant.\smallskip {\it Surjectivity}: Let
$f:\widetilde{C}\rightarrow C$ be the normalization morphism, and
${\cal L}\in \bar JC$.
Let us denote by
$\widetilde{\cal L}$ the line bundle on $\widetilde{C}$ quotient
of $f^*{\cal L}$ by its torsion subsheaf.  We
claim that its degree is
$\le 0$: we have an exact sequence
$$0\rightarrow {\cal L}\longrightarrow f_*\widetilde{\cal
L}\longrightarrow {\cal T}\rightarrow 0$$
where ${\cal T}$ is a skyscrapersheaf supported on the singular locus of $C$,
such
that $\dim{\cal T}_x\le\delta_x$ for all $x\in
C$ ([G-P], lemma 1.1); this implies  $\chi(\widetilde{\cal L})
-\chi({\cal L})\le$ $ \chi({\cal O}_{\widetilde{C}})-\chi({\cal O}_C)$,
from which the required inequality follows. Since $\widetilde{C}$ is rational,
it
follows that $\widetilde{\cal L}^{-1}$ admits a global section whose zero set
is
contained in $\Sigma $. \ind Because of the canonical isomorphisms
$$\Hom_{{\cal
O}_C}({\cal L},{\cal O}_{\widetilde{C}})\cong \Hom_{{\cal
O}_{\widetilde{C}}}(f^*{\cal L},{\cal O}_{\widetilde{C}}) \cong \Hom_{{\cal
O}_{\widetilde{C}}}(\widetilde{\cal L},{\cal O}_{\widetilde{C}})\ ,$$
we conclude that there exists a  homomorphism $i:{\cal L} \rightarrow {\cal
O}_{\widetilde{C}}$ which is bijective outside  $\Sigma$.
 Put $n_x=\dim {\cal O}_{\widetilde{C},x}/i({\cal L}_x)$ for each
$x\in\Sigma$.
  Since $$\sum_{x\in\Sigma}n_x=\dim {\cal
O}_{\widetilde{C}}/i({\cal L})=\chi( {\cal O}_{\widetilde{C}})-\chi({\cal
L})=g=\sum_{x\in\Sigma}
\delta_x\ ,$$
there exists a rational function $\varphi$ on $\widetilde{C}$ with divisor
$\sum_x(\delta_x-n_x)\,x$. Replacing ${\cal L}$ by $\varphi{\cal L}$, we may
assume
$n_x=\delta_x$ for all $x$, which means that  ${\cal L}$ belongs to the image
of
$e$. \cqfd

\medskip
\ind The variety ${\bf G}_x$ depends only on the local ring ${\cal O}$
of $C$ at $x$ (even only on its completion); we will also denote it by ${\bf
G}_{\cal O}$. Recall that ${\bf G}_{\cal O}$ parameterizes the sub-$\!{\cal
O}$\tx
modules $L$ of the normalization $\widetilde{\cal O}$ of ${\cal O}$ with $\dim
\widetilde{\cal O}/L= \dim  \widetilde{\cal O}/{\cal O}$. We put $\varepsilon
(x)=e({\bf G}_x)$ (or
$\varepsilon ({\cal O})=e({\bf G}_{\cal O})$). The above Proposition gives:

 \th Proposition
\enonce Let $C$ be a rational  unibranch curve; then $e(\bar J C)=\pprod_{x\in
C}^{}\varepsilon (x)$.~\cqfd\endth\label{prod}
\ind Of course $\varepsilon (x)$ is equal to $1$ for a smooth point, so we
could as
well consider the product over the singular locus $\Sigma $ of $C$. Note that
in
view of Prop.\ \ref{elag}, we may define $\varepsilon (x)$ for a non-unibranch
singularity by taking the product of the $\varepsilon $\tx invariants of each
branch; Prop.\ \ref{prod} remains valid.

\section{Examples}
\subsec {\it Singularities with ${\bf C}^*$\tx action}\label{C*}
\ind Assume that the local, unibranch ring ${\cal O}$ admits a ${\bf C}^*$\tx
action. This action extends to its completion, so we will assume that
${\cal O}$ is complete. The ${\bf C}^*$\tx action also extends to the
normalization
$\widetilde{\cal O}$ of ${\cal O}$, and there exists a local coordinate
$t\in\widetilde{\cal O}$ such that the line ${\bf C}t$ is preserved (this is
because  the pro-algebraic group $\Aut(\widetilde{\cal O})$ is an extension of
${\bf C}^*$ by a pro-unipotent group, hence all subgroups of
$\Aut(\widetilde{\cal
O})$ isomorphic to ${\bf C}^*$ are conjugate). It follows that the graded
subring
${\cal O}$ is associated to a semi-group $\Gamma \i{\bf N}$, i.e.\ ${\cal O}$
is
the ring ${\bf C}[[\Gamma ]]$ of formal series $\displaystyle \sum_{\gamma
\in\Gamma }a_\gamma t^\gamma $. \ind The ${\bf C}^*$\tx actions on ${\cal O}$
and
$\widetilde{\cal O}$ give rise to a ${\bf C}^*$\tx action on
${\bf G}_{\cal O}$. The fixed points of this action are the submodules of
$\widetilde{\cal O}$ which are graded, that is of the form ${\bf C}[[\Delta
]]$,
where $\Delta $ is a subset of ${\bf N}$; the condition $\dim \widetilde{\cal
O}/{\bf C}[[\Delta ]]=\dim \widetilde{\cal O}/{\cal O}$ means $\Card({\bf
N}\moins\Delta )= \Card({\bf N}\moins\Gamma) $, and the condition that ${\bf
C}[[\Delta ]]$ is a ${\cal O}$\tx module means $\Gamma +\Delta \i\Delta $. The
first condition already implies that there are only finitely many such fixed
points. According to [B], the number of these fixed points is equal to $e({\bf
G}_{\cal O})$. We conclude: \th Proposition
\enonce Let $\Gamma \i{\bf N}$ be a semi-group with finite complement. The
number
$\varepsilon ({\bf C}[[\Gamma ]])$ is equal to the number of subsets $\Delta
\i{\bf N}$ such that $\Gamma +\Delta \i\Delta $ and $\Card({\bf N}\moins\Delta
)=
\Card({\bf N}\moins\Gamma) $. \cqfd
\endth\label{\Gamma}
\ind I do not know whether there exists a closed formula computing this number,
say in terms of a minimal set of generators of $\Gamma $. This turns out to be
the
case in the situation we were originally interested in, namely planar
singularities. The semi-group $\Gamma $ is then generated by two coprime
integers
$p$ and $q$, which means that the local ring ${\cal O}$ is of the form ${\bf
C}[[u,v]]/(u^p-v^q)$.

\th Proposition
\enonce Let $p,q$ be two coprime integers. Then $$ \varepsilon ({\bf
C}[[u,v]]/(u^p-v^q))= {1\over p+q}{p+q\choose p}\ .$$
\endth\label{pq}

{\it Proof}: The following proof has been shown to me by P. Colmez.
\ind (\ref{pq}.1) We first observe that if a subset $\Delta$ satisfies $\Gamma
+\Delta \i\Delta $, all its translates $n+\Delta$ $(n\in{\bf Z})$ contained in
${\bf N}$ have the same property; moreover, among all these translates, there
is
exactly one
 with $\Card({\bf N}\moins\Delta ) = \Card({\bf N}\moins\Gamma)$. Thus the
number we want to compute is the cardinal of the set ${\cal D}$ of subsets
$\Delta\i{\bf N}$ such that $\Gamma
+\Delta \i\Delta $, modulo the identification of a subset and its translates.

\smallskip
\ind (\ref{pq}.2) For such a subset $\Delta$, let us introduce the
generating function
$\displaystyle F_\Delta (T)=  \sum_{\delta\in\Delta }T^\delta\in {\bf
Z}[[T]]$. Since
$p+\Delta\i\Delta $, we can write, in a unique way, $ \Delta
=\union_{i=1}^p (a(i)+p{\bf N})$ ; then
 $\displaystyle\ (1-T^p)\,F_\Delta (T)=\sum_{i=1}^p T^{a(i)}$ .
 Writing similarly $ \Delta =$ $\union_{j=1}^q (b(j)+q{\bf
N})$, we get $\displaystyle\ (1-T^q)\,F_\Delta (T)=\sum_{j=1}^q T^{b(j)}$ . Put
$a(j)=b(j-p)+p$ for $p+1\le j\le p+q$;
 the equality
$\displaystyle\ (1-T^p)\sum_{j=p+1}^{p+q} T^{a(j)-p}=(1-T^q)\sum_{i=1}^p
T^{a(i)}\ $ reads
$$\sum_{i=1}^{p+q} T^{a(i)} = \sum_{i=1}^p T^{a(i)+q} +
\sum_{j=p+1}^{p+q} T^{a(j)-p}\quad.\leqno (\ref{pq}\ {\it a})$$
\ind Conversely, given a function $a:[1,p+q]\rightarrow {\bf N}$ satisfying
(\ref{pq} {\it a}), the set $ \Delta =$ $\union_{i=1}^p (a(i)+p{\bf N})$ is
equal
to $\union_{j=p+1}^{p+q}(a(j)-p+q{\bf N})$, and therefore satisfies $\Gamma
+\Delta \i \Delta $ (note that (\ref{pq}~{\it a}) implies that the classes
(mod.\
$p$) of the $a(i)$'s, for $1\le i\le p$, are all distinct).
\ind  The  equality (\ref{pq} {\it a}) means that there
exists a permutation $\sigma \in{\goth S}_{p+q}$ such that $a(\sigma i)$ is
equal
to $a(i)+q$ if $i\le p$ and to $a(i)-p$ if $i>p$. This implies that $a(\sigma
^m(i))$ is of the form $a(i)+\alpha q-\beta p$, with $\alpha ,\beta \in{\bf N}$
and $\alpha +\beta =m$; since $p$ and $q$ are coprime, it follows that $\sigma
$
is of order $p+q$, i.e.\ is a circular permutation. It also follows that the
numbers $a(i)$ are all distinct; hence the permutation $\sigma $ is uniquely
determined. Let $\tau $ be a permutation such that $\tau \sigma \tau^{-1}$  is
the
permutation  $i\mapsto i+1$ (mod.\ $p+q)$, and let $S_\Delta =\tau ([1,p])$.
Replacing $a$ by $a\rond\tau^{-1} $, our
 function $a$ satisfies
$$a(i+1)=\cases{
a(i)+q  & if $i\in S_\Delta $,\cr
a(i)-p  & if $i\notin S_\Delta $\ .}\leqno(\ref{pq}\ {\it b})$$
Since $\tau $ is determined up to right multiplication by a power of $\sigma $,
the
set $S_\Delta \i[1,p+q]$ is well determined up to a translation (mod.\ $p+q)$.
Note that replacing $\Delta$ by $n+\Delta$ amounts to add the constant value
$n$ to
the function
$a$, hence does not change  $S_\Delta$.
\ind (\ref{pq}.3) Conversely, let us start from  a subset $S\i[1,p+q]$ with $p$
elements. We
define inductively a function $a_S$ on $[1,p+q]$ by the relations (\ref{pq}
{\it
b}),  giving to $a_S(1)$ an arbitrary value, large enough so that $a_S$
takes its values in ${\bf N}$.   By construction the function $a_S$ satisfies
(\ref{pq} {\it b}), so by (\ref{pq}.2) the subset
 $ \Delta_S =\union_{s\in S}^{} (a_S(s)+p{\bf N})$   satisfies $\Gamma
+\Delta_S \i\Delta_S $.
\ind An easy computation gives $a_{S+1}(i+1)=a_S(i)$ and therefore
$\Delta _{S+1}=\Delta _S$. Let ${\cal S}$ be the set of subsets of $[1,p+q]$
with
$p$ elements, modulo translation; the maps
 $\Delta \mapsto S_\Delta $ from ${\cal D}$ to ${\cal S}$ and $S\mapsto
\Delta_S $ from ${\cal S}$ to ${\cal D}$ are inverse of each other. Since
$\displaystyle \Card({\cal S})={1\over
p+q}{p+q\choose p}$, the Proposition follows. \cqfd

\bigskip
\subsec {\it Simple singularities}
\ind We now consider the case where the singularities of $C$ are {\it simple},
i.e.\
of $A,D,E$ type. The local ring of
such a singularity has only finitely many isomorphism classes of torsion free
rank
$1$ modules, and this property characterizes these singularities among all
plane
curves singularities [G-K].

\th Proposition \enonce Let ${\cal O}$ be the local ring of a
simple singularity. Then $\varepsilon ({\cal O})$ is the number of isomorphism
classes of torsion free rank $1\ {\cal O}$\tx modules. It is given by:
 \ind -- $\varepsilon ({\cal O})=l+1$ if ${\cal O}$ is of type $A_{2l}\,;$ \ind
--
$\varepsilon ({\cal O})=1$ if ${\cal O}$ is of type $A_{2l+1}\,;$ \ind --
$\varepsilon ({\cal O})=1$ if ${\cal O}$ is of type $D_{2l}\ (l\ge2)\,;$ \ind
--
$\varepsilon ({\cal O})=l$ if ${\cal O}$ is of type $D_{2l+1}\ (l\ge2)\,;$ \ind
--
$\varepsilon ({\cal O})=5$ if ${\cal O}$ is of type $E_6\,;$ \ind --
$\varepsilon
({\cal O})=2$ if ${\cal O}$ is of type $E_7\,;$
 \ind -- $\varepsilon ({\cal O})=7$ if ${\cal O}$ is of type $E_8$.
\endth\label{rat}
{\it Proof}: Let $C$ be a rational curve with only one simple singularity, with
local ring ${\cal O}$; the action of $JC$ on $\bar JC$ has finitely many
orbits,
corresponding to the different  isomorphism classes of  rank  $1\ {\cal O}$\tx
modules. Since each  orbit is an affine space, its Euler characteristic is $1$,
hence by (\ref{e}) $ \varepsilon ({\cal O})=e(\bar JC)$ is equal to the number
of
these orbits.
\ind If ${\cal O}$ is unibranch, its completion is of the form
${\bf C}[[u,v]]/(u^p-v^q)$,
with $p=2$, $q=2l+1$ for the type $A_{2l}$, $p=3$, $q=4$  for the type $E_6$
and
$p=3$, $q=5$ for the type $E_8$; in these cases the result follows from
\ref{pq}.
 We have
already observed that $\varepsilon =1$ for a $A_{2l+1}$ singularity (Remark
\ref{An}). A $D_l$ singularity is the union of a $A_{l-3}$ branch and a
transversal smooth branch, hence the result by \ref{elag}. Finally an $E_7$
singularity is the union of an ordinary cusp and its tangent, hence has
$\varepsilon =2$. \cqfd

\smallskip
\rem{Remark} Let ${\cal D}$ be the set  of graded sub-$\!{\cal O}$\tx modules
$L\i\widetilde{\cal O}$ with $\dim \widetilde{\cal O}/L=$ $\dim \widetilde{\cal
O}/{\cal O}$. Two  modules $L$, $M$ in ${\cal D}$ are isomorphic if and only if
$M=t^nL$ for some
$n\in{\bf Z}$, but the dimension condition forces $n=1$. It follows that {\it
each
torsion free rank $1\ {\cal O}$\tx module is isomorphic to exactly one element
of}
${\cal D}$.
 It is quite easy that way to write down the list of isomorphism classes of
rank $1\ {\cal O}$\tx modules (which is of course well-known, see e.g.\ [G-K]).
For
instance if ${\cal O}$ is of type $E_8$, we get the following
modules (with the notation of \ref{C*}):\par
${\cal O}$, ${\cal O}t+{\cal O}t^8$, ${\cal O}t^2+{\cal O}t^6$, ${\cal O}t^2+
{\cal O}t^4$, ${\cal O}t^3+{\cal O}t^4$, ${\cal O}t^3+{\cal O}t^5+{\cal O}t^7$,
$\widetilde{\cal O}t^4$.

\vskip2cm
\centerline{ REFERENCES} \vglue15pt\baselineskip12.8pt
\def\num#1{\smallskip \item{\hbox to\parindent{\enskip [#1]\hfill}}}
\parindent=1.3cm
\num{B} A.\ {\pc BIALYNICKI-BIRULA}: {\sl On fixed point schemes of actions of
multiplicative and additive groups}. Topology {\bf 12}, 99-103  (1973).
 \num{G}  L.\ {\pc G\"OTTSCHE}: {\sl The Betti numbers of the Hilbert scheme of
points on a smooth projective surface}. Math.\ Ann.\ {\bf 86}, 193-207 (1990).
\num{G-K} G.-M.\ {\pc GREUEL}, H.\ {\pc KN\"ORRER}: {\sl Einfache
Kurvensingularit\"aten und torsionfreie Moduln}. Math.\ Ann.\ {\bf 270},
417-425
(1985). \num{G-P} G.-M.\ {\pc GREUEL}, G.\ {\pc PFISTER}: {\sl Moduli spaces
for
torsion free modules on curve singularities, I}. J.\ Algebraic Geometry {\bf
2},
81-135 (1993).
\num{H} D. {\pc HUYBRECHTS}: {\sl Birational symplectic manifolds and their
deformations}. Preprint alg-geom/9601015.
 \num{M} S.\ {\pc MUKAI}: {\sl Symplectic structure of the moduli
space of sheaves on an abelian or K3 surface}. Invent.\ math.\ {\bf 77},
101-116
(1984).
\num{R} C. J. {\pc REGO}: {\sl The compactified Jacobian}. Ann.\ scient.\
\'Ec.\
Norm.\ Sup.\ {\bf 13}, 211-223 (1980).
  \num{V} J.-L.\ {\pc VERDIER}: {\sl Stratifications de Whitney et th\'eor\`eme
de
Bertini-Sard}. Invent.\ math.\ {\bf 36}, 295-312 (1976).
\num{Y-Z} S.-T.\ {\pc YAU}, E.\ {\pc ZASLOW}: {\sl BPS states, string duality,
and
nodal curves on K3}. Preprint hep-th/9512121.

\vskip1cm
\def\pc#1{\eightrm#1\sixrm}
\hfill\vtop{\eightrm\hbox to 5cm{\hfill Arnaud {\pc BEAUVILLE}\hfill}
 \hbox to 5cm{\hfill DMI -- \'Ecole Normale
Sup\'erieure\hfill} \hbox to 5cm{\hfill (URA 762 du CNRS)\hfill}
\hbox to 5cm{\hfill  45 rue d'Ulm\hfill}
\hbox to 5cm{\hfill F-75230 {\pc PARIS} Cedex 05\hfill}}

\end